\documentclass[a4paper, 12pt]{amsproc}
\allowdisplaybreaks[1]

\usepackage{amsmath, amsthm, amssymb, mathtools, mathrsfs, stmaryrd}
\usepackage{ascmac}
\usepackage{comment}
\usepackage{bm}
\usepackage{ytableau}

\usepackage{hyperref}

\usepackage{tikz}

\usepackage{graphicx}
\usepackage{here}
\usepackage{time}
\usepackage[abbrev]{amsrefs}

\usepackage{xcolor}
\usepackage[capitalize,nameinlink,noabbrev,nosort]{cleveref}
\hypersetup{
 colorlinks=true,       
 linkcolor=blue,          
 citecolor=blue,        
 filecolor=blue,      
 urlcolor=blue,           
}

\usepackage{fullpage}

\makeatletter
\@namedef{subjclassname@2020}{%
  \textup{2020} Mathematics Subject Classification}
\makeatother

\newtheorem{theoremcounter}{Theorem Counter}[section]

\theoremstyle{definition}

\newtheorem{remark}[theoremcounter]{Remark}
\newtheorem{example}[theoremcounter]{Example}

\theoremstyle{plain}

\newtheorem{proposition}[theoremcounter]{Proposition}
\newtheorem{corollary}[theoremcounter]{Corollary}

\newtheorem{theorem}[theoremcounter]{Theorem}

\numberwithin{equation}{section}

\newcommand{\R}{\mathbb{R}}

\DeclareMathOperator{\ReNew}{Re}
\renewcommand{\Re}{\ReNew}

\begin{document}
\author{Takashi Miyagawa}
\address[Takashi Miyagawa]{Onomichi City University,  1600-2 Hisayamada-cho, Onomichi, Hiroshima, 722-8506, Japan} 
\email{miyagawa@onomichi-u.ac.jp}



\begin{abstract}
We investigate the Laurent series expansions of the Barnes double zeta function
$\zeta_2(s,\alpha;v,w)$ at $s=1$ and $s=2$.
We derive explicit limit representations for the Laurent coefficients, which may be regarded as analogues of the Euler--Stieltjes constants in the Barnes setting. In particular, we obtain formulas in terms of finite double sums with explicit correction terms, including a new limit representation for the constant term at $s=1$. We also establish relations between the Laurent coefficients at $s=1$ and the Taylor coefficients at $s=0$, as well as relations between the coefficients at $s=1$ and $s=2$.

Furthermore, we study the asymptotic behavior of these coefficients as the order tends to infinity and show that they approach simple closed-form expressions with exponentially decaying error terms. Finally, we investigate arithmetic properties of certain constant terms in the Laurent expansions. We establish transcendence results for some infinite families of special values, including a family at $s=1$ for which every member is transcendental and a family at $s=2$ for which all but at most one are transcendental.

\end{abstract}


\title{On the Laurent series expansions of the Barnes double zeta function}

\maketitle

\section{Introduction}

Let $r$ be a positive integer and let $s=\sigma+it \ (\sigma,t\in \R)$ be a complex variable.
For $ \alpha>0$ and $w_1,\dots,w_r>0$, the Barnes multiple zeta function, introduced by E.~W.~Barnes \cites{Barnes1899, Barnes1901, Barnes1904}, is defined by
\begin{align}\label{zeta_r}
\zeta_r(s,\alpha ;w_1,\dots,w_r)
=
\sum_{m_1=0}^\infty \cdots \sum_{m_r=0}^\infty
\frac{1}{(\alpha+m_1 w_1+\cdots+m_r w_r)^s}.
\end{align}
This series is a natural generalization of the Hurwitz zeta function
\begin{align*}\zeta_H(s,a) = \sum_{m=0}^\infty \frac{1}{(m+a)^s} \qquad (\sigma>1).
\end{align*}
The series \eqref{zeta_r} converges absolutely for $\sigma>r$ and admits a meromorphic continuation to the entire complex plane, with possible simple poles at $s=1,2,\dots,r$.
In this paper, we focus on the case $r=2$. The Barnes double zeta function is defined by
\begin{equation}\label{B_2-zeta}
\zeta_2(s,\alpha;v,w)
=
\sum_{m=0}^\infty \sum_{n=0}^\infty
\frac{1}{(\alpha+mv+nw)^s},
\end{equation}
where $\alpha,v,w>0$. This series converges absolutely for $\sigma>2$ and extends to a meromorphic function on $\mathbb{C}$ with a simple pole at $s=2$ and at most a simple pole at $s=1$. The residues are given by
\begin{align*}\operatorname{Res}_{s=2}\zeta_2(s,\alpha;v,w)=\frac{1}{vw},
\qquad
\operatorname{Res}_{s=1}\zeta_2(s,\alpha;v,w)=\frac{v+w-2\alpha}{2vw},
\end{align*}
which are classical results due to Barnes \cite{Barnes1901}.

The Barnes double zeta function is closely related to the Barnes double gamma function, defined by
\[
\log \Gamma_2(x;v,w)
:=
\left.
\frac{\partial}{\partial s}
\zeta_2(s,x;v,w)
\right|_{s=0}
=
\zeta_2'(0,x;v,w).
\]
Its logarithmic derivative with respect to $x$ is called the Barnes double digamma function:
\[
\psi_2(x;v,w)
:=
\frac{d}{dx}\log\Gamma_2(x;v,w)
=
\frac{d}{dx}\zeta_2'(0,x;v,w).
\]
This function is related to the constant term in the Laurent expansion of
$\zeta_2(s,\alpha;v,w)$ at $s=1$.

To place our results in context, we recall the classical Laurent expansions of the Riemann and Hurwitz zeta functions. The Riemann zeta function
\[
\zeta(s)=\sum_{m=1}^\infty \frac{1}{m^s} \qquad (\sigma>1)
\]
admits the Laurent expansion at $s=1$
\begin{equation*}\zeta(s)=\frac{1}{s-1}+\gamma+\sum_{k=1}^\infty \gamma_k (s-1)^k,
\end{equation*}
where $\gamma$ is the Euler constant and $\gamma_k$ are the Euler--Stieltjes constants, given by
\[
\gamma_k = \frac{(-1)^k}{k!} \lim_{M\to\infty}
\left(
\sum_{m=1}^M \frac{\log^k m}{m}
-
\frac{\log^{k+1} M}{k+1}
\right).
\]
Similarly, the Hurwitz zeta function admits the expansion
\[
\zeta_H(s,a)=\frac{1}{s-1}+\gamma_0(a)+\sum_{k=1}^\infty \gamma_k(a)(s-1)^k,
\]
where
\[
\gamma_k(a)=\frac{(-1)^k}{k!} \lim_{M\to\infty}
\left(
\sum_{m=0}^M \frac{\log^k(m+a)}{m+a}
-
\frac{\log^{k+1}(M+a)}{k+1}
\right).
\]
These coefficients reflect the analytic structure of the zeta functions near their poles and play an important role in analytic number theory.
Various quantitative results, including explicit upper bounds and asymptotic estimates, are known for the Euler--Stieltjes constants and their Hurwitz analogues.
In the present paper, the symbols $\gamma_k$ and $\gamma_k(a)$ denote the coefficients of the corresponding powers in the Laurent expansions, rather than the classical Stieltjes constants before division by $k!$. With this normalization, the estimate of Zhang and Williams \cite{ZhangWilliams1994} for the classical Euler--Stieltjes constants is equivalent to
\begin{equation}
    |\gamma_k|
    \le
    \frac{3+(-1)^k}{(2\pi)^k}
    \cdot
    \frac{(2k)!}{k!\,k^{k+1}}
    \qquad (k\ge 1).
    \label{estim_gamma_k}
\end{equation}
For a discussion of equivalent estimates for the Stieltjes constants,
see also Finch \cite{Finch2003}*{\S 2.21}.
Similarly, rewriting the estimate of Berndt \cite{Berndt1972} according to the normalization adopted here, we have
\begin{equation}
    \left|
        \gamma_k(a)
        -
        \frac{(-1)^k(\log a)^k}{a k!}
    \right|
    \le
    \frac{3+(-1)^k}{k\pi^k}
    \qquad
    (0<a\le 1,\ k\ge 1).
    \label{estim_gamma_k(a)}
\end{equation}
These estimates illustrate that the Laurent coefficients contain delicate information concerning the growth of zeta functions and the analytic behavior near their singularities.

The aim of this paper is to investigate the Laurent expansions of the Barnes double zeta function at $s=1$ and $s=2$, and to obtain explicit expressions for the coefficients. These coefficients can be regarded as analogues of the Euler--Stieltjes constants in the Barnes setting. We establish explicit limit representations and clarify their structural relations with derivatives of the zeta function at $s=0$. We also investigate arithmetic properties of some special values and obtain transcendence results for infinite families of constant terms at $s=1$ and $s=2$: every member of the former family is transcendental, while all but at most one member of the latter family are transcendental.
Furthermore, we show that, in contrast to the classical Hurwitz case, these coefficients exhibit a remarkably simple asymptotic behavior. This behavior is governed by the contribution of the other pole to each Laurent expansion.

As related work, Matsumoto, Onozuka, and Wakabayashi \cite{MatsumotoOnozukaWakabayashi2020} studied multivariable Laurent-type expansions of the Euler--Zagier multiple zeta function
\[
\zeta_r(s_1, \dots, s_r)
=
\sum_{m_1=1}^\infty \cdots \sum_{m_r=1}^\infty
\frac{1}{
m_1^{s_1}
(m_1+m_2)^{s_2}
\cdots
(m_1+\cdots+m_r)^{s_r}
},
\]
and introduced analogues of Euler--Stieltjes constants in that setting. Their results provide a complementary viewpoint to the present study.
Sakane and Aoki \cite{SakaneAoki2022} obtained a reduction formula for the Barnes multiple zeta function with positive rational periods, expressing it, up to an explicit elementary factor, as a finite linear combination of Hurwitz zeta functions. We use its double-zeta specialization in Section~6.

The paper is organized as follows. In Section~2, we establish Euler--Stieltjes-type
limit formulas for the Laurent coefficients at $s=2$ and for the constant term at
$s=1$. Section~3 describes relations among the Laurent coefficients and the Taylor
coefficients at $s=0$. In Section~4, we study the asymptotic behavior of the Laurent
coefficients as $k\to\infty$. Section~5 gives a Mellin-type integral representation for the Laurent
coefficients at $s=2$.
Finally, in Section~6, we record some simple special values and their arithmetic
consequences.

\section{Euler--Stieltjes-type limit formulas}

We describe our results on the principal parts and constant terms in the Laurent expansions of the Barnes double zeta function (\ref{B_2-zeta}).
When the Laurent expansions of the Barnes double zeta function at $s=2$ and $s=1$ are given by
\begin{align}
& \zeta_2(s, \alpha; v, w) = \frac{1}{vw}\frac{1}{s-2}   + \gamma_0(2,\alpha;v,w) + \sum_{k=1}^\infty \gamma_k(2,\alpha;v,w)(s-2)^k,\label{Laurent_s=2}\\
&\zeta_2(s, \alpha; v, w)
= \frac{v+w-2\alpha}{2vw} \frac{1}{s-1}   + \gamma_0(1,\alpha;v,w) + \sum_{k=1}^\infty \gamma_k(1,\alpha;v,w)(s-1)^k.
\label{Laurent_s=1}
\end{align}
we obtain the following Euler--Stieltjes-type limit formulas.

\subsection{The Laurent expansion at $s=2$}
\begin{theorem}\label{th:Main_Theorem1}
    Let $\alpha,v,w>0$. Then, for each $k=0,1,2,\ldots$, the Laurent coefficient
    $\gamma_k(2,\alpha;v,w)$ is given by
    \begin{align}
    &\gamma_k(2,\alpha;v,w)
    =
    \lim_{M\to\infty}
    \frac{(-1)^k}{k!}
    \Bigg(
    \sum_{m=0}^M\sum_{n=0}^M
    \frac{\log^k(\alpha+mv+nw)}
         {(\alpha+mv+nw)^2}
    \nonumber\\
    &\qquad
    -
    \frac{k!}{vw}
    \bigg\{
    1+
    \sum_{j=1}^{k+1}
    \frac{1}{j!}
    \left(
    \log^j(\alpha+vM)
    +
    \log^j(\alpha+wM)
    -
    \log^j(\alpha+vM+wM)
    \right)
    \bigg\}
    \Bigg).
    \label{Euler_Stieltjes_2}
    \end{align}
\end{theorem}

\begin{example}
When $k=0$ in \eqref{Euler_Stieltjes_2}, the constant term is given by
\begin{align*}
\gamma_0(2,\alpha;v,w)
=
\lim_{M\to\infty}
\Bigg\{
&\sum_{m=0}^{M}\sum_{n=0}^{M}
\frac{1}{(\alpha+mv+nw)^2}
-
\frac{1}{vw}
\log{\frac{e(\alpha+vM)(\alpha+wM)}{(\alpha+vM+wM)}}
\Bigg\}.
\end{align*}
This formula is analogous to the classical limit representations of
the Euler constant and the generalized Euler constant.
\end{example}

\begin{proof}[Proof of Theorem \ref{th:Main_Theorem1}]
The idea of the proof is to approximate the double series by a finite sum and apply the Euler-Maclaurin summation formula. The coefficients $ \gamma_{k-1}(2, \alpha;v,w) \ (k=1,2,3,\cdots)$ of the Laurent series expansions of $ \zeta_2(s, \alpha;v,w) $ 
	are given by
	\[
	\gamma_{k-1}(2, \alpha;v,w) = \frac{1}{k!} 
	\left[\frac{d^k}{ds^k} \left\{ (s-2) \zeta_2(s,\alpha;v,w) \right\} \right]_{s=2}.
	\]
	Therefore, instead of the double infinite series $ \zeta_2(s, \alpha;v,w) $, we consider the finite partial sum
	\[
	\tilde{\zeta}_2(s,\alpha;v,w)
	= \sum_{m=0}^M \sum_{n=0}^M \frac{1}{(\alpha +mv + nw)^s}
	\]
	up to the $M$-th term, we compute the higher-order derivatives as follows:
	\begin{align*}
		& (-1)^k \frac{d^k}{ds^k}
		\left\{ (s-2) \tilde{\zeta}_2(s,\alpha;v,w) \right\} \\
		& \qquad = (s-2) \sum_{m=0}^M \sum_{n=0}^M \frac{ \log^k(\alpha+mv+nw)}{(\alpha +mv + nw)^s}
		- k \sum_{m=0}^M \sum_{n=0}^M \frac{ \log^{k-1}(\alpha+mv+nw)}{(\alpha +mv + nw)^s}.
	\end{align*}
Assume $\sigma>2$. Recall that, for a positive integer $M$ and
$f\in C^1([0,M])$, the first-order Euler--Maclaurin summation formula gives
\begin{equation}
\sum_{n=1}^M f(n)
=
\int_0^M f(x)\,dx
+
\int_0^M \{x\}f'(x)\,dx.
\label{eq:EM}
\end{equation}
Here, $\{x\}$ denotes the fractional part of $x$, that is, $\{x\}=x-[x]$.
We apply \eqref{eq:EM} to
$f(x,y)=\log^k(\alpha+vx+wy)/(\alpha+vx+wy)^s$
in each variable.
Let $M$ be a sufficiently large positive integer and consider the finite double sum
\[
\sum_{m=0}^M \sum_{n=0}^M \frac{\log^k(\alpha+mv+nw)}{(\alpha+mv+nw)^s}.
\]
In the following Euler--Maclaurin expansions, terms with vanishing prefactors
are understood to be omitted when $k=0$ or $k=1$.
First, we consider the single sum, so we have
\begin{align*}
 & \sum_{m=0}^M \frac{\log^k(\alpha+mv+nw)}{(\alpha+mv+nw)^s} \\
 & = \frac{\log^k(\alpha+nw)}{(\alpha+nw)^s}
  + \int_0^M\frac{\log^k(\alpha+vx+nw)}{(\alpha+vx+nw)^s}dx \\
 & \quad +\int_0^M\frac{\partial}{\partial x}\left\{\frac{\log^k(\alpha+vx+nw)}{(\alpha+vx+nw)^s}\right\}\{x\}dx \\
 & =\frac{\log^k(\alpha+nw)}{(\alpha+nw)^s}
  + 
  \int_0^M\frac{\log^k(\alpha+vx+nw)}{(\alpha+vx+nw)^s}dx \\
 & \quad +  kv\int_0^M\frac{\log^{k-1}(\alpha+vx+nw)}{(\alpha+vx+nw)^{s+1}}\{x\}dx
-
vs\int_0^M\frac{\log^{k}(\alpha+vx+nw)}{(\alpha+vx+nw)^{s+1}}\{x\}dx.
\end{align*}
Then we have
\begin{align}
& \sum_{m=0}^M \sum_{n=0}^M \frac{\log^k(\alpha+mv+nw)}{(\alpha+mv+nw)^s} \nonumber\\
&=\frac{\log^k\alpha}{\alpha^s}+ \int_0^M\frac{\log^k(\alpha+wy)}{(\alpha+wy)^s}dy \nonumber\\
& \quad +kw\int_0^M\frac{\log^{k-1}(\alpha+wy)}{(\alpha+wy)^{s+1}}\{y\}dy
  -sw\int_0^M\frac{\log^{k}(\alpha+wy)}{(\alpha+wy)^{s+1}}\{y\}dy\nonumber\\
& \quad 
  + \int_0^M\frac{\log^k(\alpha+vx)}{(\alpha+vx)^s}dx 
  + \int_0^M\int_0^M\frac{\log^{k}(\alpha+vx+wy)}{(\alpha+vx+wy)^{s}}dxdy \nonumber\\
& \quad 
  + kw\int_0^M\int_0^M \frac{\log^{k-1}(\alpha+vx+wy)}{(\alpha+vx+wy)^{s+1}}\{y\}dxdy \nonumber\\
& \quad
  - s w \int_0^M\int_0^M\frac{\log^{k}(\alpha+vx+wy)}{(\alpha+vx+wy)^{s+1}}\{y\}dxdy \nonumber\\
& \quad
  + kv\int_0^M \frac{\log^{k-1}(\alpha+vx)}{(\alpha+vx)^{s+1}}\{x\}dx 
  + kv\int_0^M\int_0^M\frac{\log^{k-1}(\alpha+vx+wy)}{(\alpha+vx+wy)^{s+1}}\{x\}dxdy \nonumber\\
& \quad
  + k(k-1)vw \int_0^M\int_0^M \frac{\log^{k-2}(\alpha+vx+wy)}{(\alpha+vx+wy)^{s+2}}\{x\}\{y\} dxdy \nonumber\\
& \quad
  - kvw(s+1) \int_0^M\int_0^M\frac{\log^{k-1}(\alpha+vx+wy)}{(\alpha+vx+wy)^{s+2}}\{x\}\{y\} dxdy \nonumber\\
& \quad
  - vs \int_0^M \frac{\log^{k}(\alpha+vx)}{(\alpha+vx)^{s+1}}\{x\}dx
  - vs \int_0^M\int_0^M\frac{\log^{k}(\alpha+vx+wy)}{(\alpha+vx+wy)^{s+1}}\{x\}dxdy \nonumber\\
& \quad
  - kvws \int_0^M\int_0^M \frac{\log^{k-1}(\alpha+vx+wy)}{(\alpha+vx+wy)^{s+2}}\{x\}\{y\} dxdy \nonumber\\
& \quad
  + vws(s+1) \int_0^M\int_0^M\frac{\log^{k}(\alpha+vx+wy)}{(\alpha+vx+wy)^{s+2}}\{x\}\{y\}dxdy.
  \label{sum_int}
\end{align}
\noindent
We now derive the limit representation of the Laurent coefficients.
For $\Re s>2$, put
\begin{align*}
    J(s,\alpha;v,w)
    \coloneqq
    \int_0^\infty\int_0^\infty
    \frac{1}{(\alpha+vx+wy)^s}\,dxdy.
\end{align*}
A direct calculation gives
\begin{align}
    J(s,\alpha;v,w)
    =
    \frac{\alpha^{2-s}}
    {vw(s-1)(s-2)}.
    \label{eq:J}
\end{align}
We also put
\begin{align*}
    F(s,\alpha;v,w)
    \coloneqq
    \zeta_2(s,\alpha;v,w)-J(s,\alpha;v,w).
\end{align*}

We first show that $F(s,\alpha;v,w)$ is holomorphic in a
neighborhood of $s=2$. For $k\geq 0$, equation \eqref{sum_int}
implies that
\begin{align}
&\sum_{m=0}^M\sum_{n=0}^M
\frac{\log^k(\alpha+mv+nw)}
     {(\alpha+mv+nw)^s}
-
\int_0^M\int_0^M
\frac{\log^k(\alpha+vx+wy)}
     {(\alpha+vx+wy)^s}\,dxdy
\label{eq:finite-difference}
\end{align}
is equal to the sum of all terms on the right-hand side of
\eqref{sum_int} except for the principal double integral.
Let $U$ be a sufficiently small neighborhood of $s=2$ contained
in the half-plane $\Re{(s)}>1$. Since
$ 0\leq \{x\}<1 $ and $0\leq \{y\}<1$,
the integrals occurring in the remaining terms of \eqref{sum_int}
converge uniformly on compact subsets of $U$ as $M\to\infty$.
Consequently, the expression in \eqref{eq:finite-difference}
converges locally uniformly in $U$ to a function holomorphic at
$s=2$.

On the other hand, for $\Re(s)>2$, termwise differentiation gives
\begin{align*}
&\lim_{M\to\infty}
\left\{
\sum_{m=0}^M\sum_{n=0}^M
\frac{\log^k(\alpha+mv+nw)}
     {(\alpha+mv+nw)^s}
-
\int_0^M\int_0^M
\frac{\log^k(\alpha+vx+wy)}
     {(\alpha+vx+wy)^s}\,dxdy
\right\}
\\
&\qquad=
(-1)^k
\frac{d^k}{ds^k}
F(s,\alpha;v,w).
\end{align*}
It follows by analytic continuation that $F(s,\alpha;v,w)$ is
holomorphic in a neighborhood of $s=2$, and
\begin{align}
(-1)^kF^{(k)}(2,\alpha;v,w)
&=
\lim_{M\to\infty}
\left\{
\sum_{m=0}^M\sum_{n=0}^M
\frac{\log^k(\alpha+mv+nw)}
     {(\alpha+mv+nw)^2}
\right.
\nonumber\\
&\hspace{42mm}\left.
-
\int_0^M\int_0^M
\frac{\log^k(\alpha+vx+wy)}
     {(\alpha+vx+wy)^2}\,dxdy
\right\}.
\label{eq:F-limit}
\end{align}
We next calculate the contribution of $J(s,\alpha;v,w)$ to the
Laurent expansion at $s=2$. Writing $t=s-2$, equation \eqref{eq:J}
gives
\begin{align*}
J(2+t,\alpha;v,w)
&=
\frac{1}{vw\,t}
\frac{e^{-t\log\alpha}}{1+t}
=
\frac{1}{vw\,t}
\bigg\{\sum_{p=0}^{\infty}(-1)^p t^p \bigg\}
\bigg\{\sum_{q=0}^{\infty}\frac{(-1)^q \log^{q}{\alpha}}{q!} t^q \bigg\}\\
&=
\frac{1}{vw\,t}
\sum_{n=0}^\infty
(-1)^n
\left(
\sum_{j=0}^n
\frac{\log^j\alpha}{j!}
\right)t^n .
\end{align*}
Therefore, the coefficient of $(s-2)^k$ in the Laurent expansion
of $J(s,\alpha;v,w)$ at $s=2$ is
\begin{align}
\frac{(-1)^{k+1}}{vw}
\sum_{j=0}^{k+1}
\frac{\log^j\alpha}{j!}.
\label{eq:J-coefficient}
\end{align}
Since
$
\zeta_2(s,\alpha;v,w)
=
J(s,\alpha;v,w)+F(s,\alpha;v,w)
$,
we obtain from \eqref{eq:F-limit} and \eqref{eq:J-coefficient}
that
\begin{align}
\gamma_k(2,\alpha;v,w)
&=
\frac{(-1)^k}{k!}
\lim_{M\to\infty}
\left\{
\sum_{m=0}^M\sum_{n=0}^M
\frac{\log^k(\alpha+mv+nw)}
     {(\alpha+mv+nw)^2}
\right.
\nonumber\\
&\hspace{38mm}\left.
-
\int_0^M\int_0^M
\frac{\log^k(\alpha+vx+wy)}
     {(\alpha+vx+wy)^2}\,dxdy
\right\}
\nonumber\\
&\quad+
\frac{(-1)^{k+1}}{vw}
\sum_{j=0}^{k+1}
\frac{\log^j\alpha}{j!}.
\label{eq:gamma-integral}
\end{align}
It remains to evaluate the finite double integral in
\eqref{eq:gamma-integral}. For $\Re s>2$, direct integration gives
\begin{align}
&\int_0^M\int_0^M
\frac{1}{(\alpha+vx+wy)^s}\,dxdy
\nonumber\\
&\quad=
\frac{1}{vw(s-1)(s-2)}
\left\{
\frac{1}{\alpha^{s-2}}
-
\frac{1}{(\alpha+vM)^{s-2}}
-
\frac{1}{(\alpha+wM)^{s-2}}
\right.
\nonumber\\
&\hspace{62mm}\left.
+
\frac{1}{(\alpha+vM+wM)^{s-2}}
\right\}.
\label{eq:finite-integral}
\end{align}
The apparent singularity at $s=2$ on the right-hand side of
\eqref{eq:finite-integral} is removable, since the expression in
braces vanishes at $s=2$.

Differentiating \eqref{eq:finite-integral} $k$ times with respect
to $s$ and then putting $s=2$, we obtain
\begin{align}
&\int_0^M\int_0^M
\frac{\log^k(\alpha+vx+wy)}
     {(\alpha+vx+wy)^2}\,dxdy
\nonumber\\
&\quad=
\frac{k!}{vw}
\sum_{j=1}^{k+1}
\frac{1}{j!}
\left\{
\log^j(\alpha+vM)
+
\log^j(\alpha+wM)
-
\log^j\alpha
\right.
\nonumber\\
&\hspace{65mm}\left.
-
\log^j(\alpha+vM+wM)
\right\}.
\label{eq:finite-log-integral}
\end{align}

Substituting \eqref{eq:finite-log-integral} into
\eqref{eq:gamma-integral}, the terms involving
$\log^j\alpha$ cancel. Hence,
\begin{align}
\gamma_k(2,\alpha;v,w)
&=
\frac{(-1)^k}{k!}
\lim_{M\to\infty}
\Bigg(
\sum_{m=0}^M\sum_{n=0}^M
\frac{\log^k(\alpha+mv+nw)}
     {(\alpha+mv+nw)^2}
\nonumber\\
&\qquad
-
\frac{k!}{vw}
\bigg\{
1+
\sum_{j=1}^{k+1}
\frac{1}{j!}
\left(
\log^j(\alpha+vM)
+
\log^j(\alpha+wM)
\right.
\nonumber\\
&\hspace{69mm}\left.
-
\log^j(\alpha+vM+wM)
\right)
\bigg\}
\Bigg).
\label{eq:gamma-final}
\end{align}

Finally, the existence and finiteness of the limit in
\eqref{eq:gamma-final} follow from \eqref{eq:F-limit}. Indeed,
equation \eqref{eq:finite-log-integral} shows that the subtracted
terms in \eqref{eq:gamma-final} contain exactly the divergent part
of the finite double sum as $M\to\infty$. The remaining part
converges to the value determined by the holomorphic function
$F(s,\alpha;v,w)$ at $s=2$. This completes the proof of
Theorem \ref{th:Main_Theorem1}.
\end{proof}

\subsection{The Laurent expansion at $s=1$}

The following theorem gives a finite double-sum representation for the constant term at $s=1$.
\begin{theorem}\label{th:limit_s1}
Let $\alpha,v,w>0$.
Then the constant term $\gamma_0(1,\alpha;v,w)$ in the Laurent
expansion of $\zeta_2(s,\alpha;v,w)$ at $s=1$ satisfies
\begin{align}
\gamma_0(1,\alpha;v,w)
&=
\lim_{M\to\infty}
\Bigg\{
\sum_{m=0}^{M}\sum_{n=0}^{M}
\frac{1}{\alpha+vm+wn} \nonumber\\
&\hspace{25mm} -
\frac{(v+w)\log(v+w)-v\log v-w\log w}{vw}(M+1)
\nonumber\\
&\hspace{30mm}
-
\frac{v+w-2\alpha}{2vw}
\log\left(
\frac{vw(M+1)}{v+w}
\right)
\Bigg\}.
\label{eq:limit_s1}
\end{align}
\end{theorem}

\begin{proof}
Put
\[
S_M(\alpha)
\coloneqq
\sum_{m=0}^{M}\sum_{n=0}^{M}
\frac{1}{\alpha+vm+wn},
\]
and let
$
B(x)\coloneqq \{x\}-1/2.
$
Adding $f(0)$ to both sides of \eqref{eq:EM} and using
$\{x\}=B(x)+1/2$, we obtain the centered form
\begin{align*}
\sum_{n=0}^{M}f(n)
=
\int_0^M f(x)\,dx
+
\frac{f(0)+f(M)}{2}
+
\int_0^M B(x)f'(x)\,dx.
\end{align*}
Applying this identity successively to
$
f(x,y)=1/(\alpha+vx+wy)
$
in the $x$- and $y$-variables, we obtain
\begin{align}
S_M(\alpha)=P_M(\alpha;v,w)+R_M(\alpha;v,w),
\label{eq:S=P+R}
\end{align}
where
\begin{align}
P_M(\alpha;v,w)
&=
\int_0^M\int_0^M
\frac{dx\,dy}{\alpha+vx+wy}
\nonumber\\
&\quad
+
\frac{1}{2w}
\log
\frac{(\alpha+wM)(\alpha+(v+w)M)}
     {\alpha(\alpha+vM)}
\nonumber\\
&\quad
+
\frac{1}{2v}
\log
\frac{(\alpha+vM)(\alpha+(v+w)M)}
     {\alpha(\alpha+wM)}
\nonumber\\
&\quad
+
\frac14
\left\{
\frac1\alpha
+\frac1{\alpha+vM}
+\frac1{\alpha+wM}
+\frac1{\alpha+(v+w)M}
\right\},
\label{eq:P_M}
\end{align}
and
\begin{align}
R_M(\alpha;v,w)
&=
-\frac{v}{w}
\int_0^M
B(x)
\left\{
\frac1{\alpha+vx}
-\frac1{\alpha+vx+wM}
\right\}dx
\nonumber\\
&\quad
-\frac{w}{v}
\int_0^M
B(y)
\left\{
\frac1{\alpha+wy}
-\frac1{\alpha+vM+wy}
\right\}dy
\nonumber\\
&\quad
-\frac{v}{2}
\int_0^M
B(x)
\left\{
\frac1{(\alpha+vx)^2}
+
\frac1{(\alpha+vx+wM)^2}
\right\}dx
\nonumber\\
&\quad
-\frac{w}{2}
\int_0^M
B(y)
\left\{
\frac1{(\alpha+wy)^2}
+
\frac1{(\alpha+vM+wy)^2}
\right\}dy
\nonumber\\
&\quad
+
2vw
\int_0^M\int_0^M
\frac{B(x)B(y)}
     {(\alpha+vx+wy)^3}\,dxdy.
\label{eq:R_M}
\end{align}

Since $B(x)$ is bounded, periodic of period one, and has mean value zero,
the first two integrals in \eqref{eq:R_M} converge by Dirichlet's test.
The remaining integrals converge absolutely. Hence
\[
R_M(\alpha;v,w)\longrightarrow R(\alpha;v,w)
\qquad (M\to\infty),
\]
where
\begin{align*}
R(\alpha;v,w)
&=
-\frac{v}{w}
\int_0^\infty
\frac{B(x)}{\alpha+vx}\,dx
-\frac{w}{v}
\int_0^\infty
\frac{B(y)}{\alpha+wy}\,dy
\\
&\quad
-\frac{v}{2}
\int_0^\infty
\frac{B(x)}{(\alpha+vx)^2}\,dx
-\frac{w}{2}
\int_0^\infty
\frac{B(y)}{(\alpha+wy)^2}\,dy
\\
&\quad
+
2vw
\int_0^\infty\int_0^\infty
\frac{B(x)B(y)}
     {(\alpha+vx+wy)^3}\,dxdy.
\end{align*}
We next relate $R(\alpha;v,w)$ to the constant term in the
Laurent expansion at $s=1$. For $\Re{(s)}>2$, applying the centered
form of \eqref{eq:EM} in both variables and letting $M\to\infty$ gives
\begin{align*}
\zeta_2(s,\alpha;v,w)
&=
\frac{\alpha^{2-s}}
     {vw(s-1)(s-2)}
+
\frac{\alpha^{1-s}}{2v(s-1)}
+
\frac{\alpha^{1-s}}{2w(s-1)}
+
\frac14\alpha^{-s}
+
\mathcal{R}(s),
\end{align*}
where
\begin{align*}
\mathcal{R}(s)
&=
-\frac{v}{w}
\int_0^\infty
\frac{B(x)}{(\alpha+vx)^s}\,dx
-\frac{w}{v}
\int_0^\infty
\frac{B(y)}{(\alpha+wy)^s}\,dy
\\
&\quad
-\frac{sv}{2}
\int_0^\infty
\frac{B(x)}{(\alpha+vx)^{s+1}}\,dx
-\frac{sw}{2}
\int_0^\infty
\frac{B(y)}{(\alpha+wy)^{s+1}}\,dy
\\
&\quad
+s(s+1)vw
\int_0^\infty\int_0^\infty
\frac{B(x)B(y)}
     {(\alpha+vx+wy)^{s+2}}\,dxdy.
\end{align*}
The function $\mathcal{R}(s)$ is holomorphic in a neighborhood of
$s=1$, and
\[
\mathcal{R}(1)=R(\alpha;v,w).
\]
Writing $s=1+h$, we have
\begin{align*}
\frac{\alpha^{2-s}}
     {vw(s-1)(s-2)}
&=
-\frac{\alpha}{vw}\frac1h
+
\frac{\alpha}{vw}(\log\alpha-1)
+O(h),
\\
\frac{\alpha^{1-s}}{2v(s-1)}
&=
\frac{1}{2v}\frac1h
-\frac{\log\alpha}{2v}
+O(h),
\\
\frac{\alpha^{1-s}}{2w(s-1)}
&=
\frac{1}{2w}\frac1h
-\frac{\log\alpha}{2w}
+O(h).
\end{align*}
Therefore,
\begin{align}
& 
\gamma_{-1}(1,\alpha;v,w)
=
\frac{v+w-2\alpha}{2vw},
\nonumber
\\
&
\gamma_0(1,\alpha;v,w)
=
\frac{\alpha}{vw}(\log\alpha-1)
-\frac{v+w}{2vw}\log\alpha
+\frac1{4\alpha}
+
R(\alpha;v,w),
\label{eq:gamma0=C0+R}
\end{align}
For convenience, put
\[
A(v,w)
:=
\frac{(v+w)\log(v+w)-v\log v-w\log w}{vw}.
\]

It remains to determine the asymptotic behavior of $P_M$.
The double integral in \eqref{eq:P_M} is
\begin{align*}
\int_0^M\int_0^M
\frac{dx\,dy}{\alpha+vx+wy}
&=
\frac1{vw}
\Bigg\{
(\alpha+(v+w)M)\log(\alpha+(v+w)M)
\\
&\qquad
-(\alpha+vM)\log(\alpha+vM)
\\
&\qquad
-(\alpha+wM)\log(\alpha+wM)
+\alpha\log\alpha
\Bigg\}.
\end{align*}
For fixed $c>0$, using
\[
(\alpha+cM)\log(\alpha+cM)
=
cM(\log M+\log c)
+
\alpha(\log M+\log c+1)
+
o(1)
\qquad (M\to\infty),
\]
we obtain from \eqref{eq:P_M}
\begin{align}
P_M(\alpha;v,w)
&=
A(v,w)M
+
\gamma_{-1}(1,\alpha;v,w)\log M
\nonumber\\
&\quad
+
A(v,w)
+
\gamma_{-1}(1,\alpha;v,w)
\log\frac{vw}{v+w}
\nonumber\\
&\quad+
\frac{\alpha}{vw}(\log\alpha-1)
-\frac{v+w}{2vw}\log\alpha
+\frac1{4\alpha}
+
o(1).
\label{eq:P-asymp}
\end{align}
Combining \eqref{eq:S=P+R}, \eqref{eq:gamma0=C0+R},
and \eqref{eq:P-asymp}, we obtain
\begin{align*}
S_M(\alpha)
&=
A(v,w)M
+
\gamma_{-1}(1,\alpha;v,w)\log M
+
A(v,w)
\\
&\quad
+
\gamma_{-1}(1,\alpha;v,w)
\log\frac{vw}{v+w}
+
\gamma_0(1,\alpha;v,w)
+
o(1).
\end{align*}
Since
\[
\log(M+1)=\log M+o(1),
\]
this can equivalently be written as
\begin{align}
S_M(\alpha)
&=
A(v,w)(M+1)
+
\gamma_{-1}(1,\alpha;v,w)
\log\left(
\frac{evw(M+1)}{v+w}
\right)
\nonumber\\
&\quad
+
\gamma_0(1,\alpha;v,w)
-
\gamma_{-1}(1,\alpha;v,w)
+
o(1).
\label{eq:S-asymp2}
\end{align}
Using $\log e=1$, rearranging \eqref{eq:S-asymp2} gives \eqref{eq:limit_s1}.
\end{proof}

\section{Relations among the Laurent coefficients}

We next describe relations among the Laurent coefficients at $s=1$ and $s=2$ and the Taylor coefficients at $s=0$.
For convenience, we write
\[
\gamma_{-1}(1,\alpha;v,w)
:=
\frac{v+w-2\alpha}{2vw}.
\]
\begin{theorem}\label{th:Main_Theorem2}
    The following relation holds between the Taylor coefficients at $s=0$ and the Laurent coefficients at $s=1$:
    \begin{align}
        \gamma_{k}(1,\alpha;v,w)
        = -\frac{1}{(k+1)!}\frac{\partial}{\partial \alpha}\zeta_2^{(k+1)}(0,\alpha;v,w) 
        \qquad (k=-1,0,1,2, \cdots) \label{th:Main_Theorem2-1}
    \end{align}
    where $ \zeta_2'(0,\alpha;v,w) = \log{\Gamma_2(\alpha;v,w)} $.
Furthermore, the following relation holds between $ \gamma_k(1,\alpha;v,w) $ and $ \gamma_k(2,\alpha;v,w) $:
\begin{align}
    \sum_{l=-1}^k (-1)^{k-l+1} \frac{\partial}{\partial \alpha}\gamma_l(1,\alpha;v,w) = \gamma_k(2,\alpha;v,w)
    \qquad (k=0,1,2,\ldots).
    \label{th:Main_Theorem2-2}
\end{align} 
\end{theorem}

\begin{example}
When $k=0, 1$ in \eqref{th:Main_Theorem2-1}, we have
\begin{align*}
    & \gamma_0(1,\alpha;v,w)=-\frac{\partial}{\partial \alpha}\zeta_2'(0,\alpha;v,w) =- \psi_2(\alpha;v,w), \\
    & \gamma_1(1,\alpha;v,w)=-\frac{1}{2}\frac{\partial}{\partial \alpha}\zeta_2''(0,\alpha;v,w).
\end{align*}
\end{example}

\begin{proof}[Proof of Theorem \ref{th:Main_Theorem2}]
Since $\zeta_2(s,\alpha;v,w)$ is holomorphic at $s=0$ and in its neighborhood, its Taylor expansion at $s=0$ is given by
\begin{align}
    \zeta_2(s,\alpha;v,w) = \sum_{k=0}^\infty \frac1{k!} \zeta_2^{(k)}(0,\alpha;v,w) s^k.
    \label{Taylor}
\end{align}
Replacing $s$ by $s+1$ in \eqref{Laurent_s=1}, we obtain
\begin{align*}
    \zeta_2(s+1, \alpha; v, w)
	= \frac{v+w-2\alpha}{2vw} \frac{1}{s}   + \gamma_0(1,\alpha;v,w) + \sum_{k=1}^\infty \gamma_k(1,\alpha;v,w)s^k.
\end{align*}
For $\Re s>2$, the defining double series and its $\alpha$-derivative
converge locally uniformly for $\alpha>0$. Termwise differentiation therefore gives
\[
\frac{\partial}{\partial\alpha}\zeta_2(s,\alpha;v,w)
=
-s\zeta_2(s+1,\alpha;v,w).
\]
The meromorphic continuation depends locally holomorphically on $\alpha$
for $\alpha>0$, as follows from the regularized Mellin representation.
Thus this identity extends meromorphically in $s$, in particular to a
neighborhood of $s=0$. The factor $s$ cancels the possible simple pole of
$\zeta_2(s+1,\alpha;v,w)$ at $s=0$, so the right-hand side is holomorphic there.
Consequently, differentiating the Taylor expansion \eqref{Taylor} with
respect to $\alpha$, we obtain
\begin{align}
     (-s) \zeta_2(s+1, \alpha; v, w) = \sum_{k=0}^\infty \frac{1}{k!}  \frac{\partial}{\partial \alpha}\zeta_2^{(k)}(0,\alpha;v,w) s^k.
     \label{Taylor_s=0}
\end{align}
On the other hand, by using \eqref{Laurent_s=1}, we have
\begin{align}
     (-s) \zeta_2(s+1, \alpha; v, w) 
     &= -s \left\{ \frac{\gamma_{-1}(1,\alpha;v,w)}{s} + \gamma_0(1,\alpha;v,w) + \sum_{k=1}^\infty \gamma_k(1,\alpha;v,w)s^k \right\} \nonumber\\
     &= - \gamma_{-1}(1,\alpha;v,w) - \gamma_0(1,\alpha;v,w)s - \sum_{k=1}^\infty \gamma_k(1,\alpha;v,w)s^{k+1}.
     \label{Laurent_s=0}
\end{align}
By the uniqueness of power series expansions, comparing the coefficients on the right-hand sides of \eqref{Taylor_s=0} and \eqref{Laurent_s=0}, we obtain
\begin{align*}
    \gamma_{k-1}(1,\alpha;v,w) = -\frac{1}{k!} \frac{\partial}{\partial \alpha} \zeta_2^{(k)}(0,\alpha;v,w)
    \quad (k=0,1,2, \dots).
\end{align*}
Replacing $k$ by $k+1$ gives \eqref{th:Main_Theorem2-1}.
Furthermore, by differentiating both sides of \eqref{Laurent_s=1} with respect to $\alpha$, we obtain
\begin{align*}
     (-s) \zeta_2(s+1, \alpha; v, w) &= \frac{ \frac{\partial}{\partial \alpha}\gamma_{-1}(1,\alpha;v,w)}{s-1} +  \frac{\partial}{\partial \alpha}\gamma_0(1,\alpha;v,w)
    \\
    &\quad+ \sum_{k=1}^\infty  \frac{\partial}{\partial \alpha}\gamma_k(1,\alpha;v,w)(s-1)^k.
\end{align*}
Changing $s$ in the above equation with $s-1$, we have
\begin{align*}
    -(s-1) \zeta_2(s, \alpha; v, w) &= \frac{ \frac{\partial}{\partial \alpha}\gamma_{-1}(1,\alpha;v,w)}{s-2} +  \frac{\partial}{\partial \alpha}\gamma_0(1,\alpha;v,w)
    \\
    &\quad+ \sum_{k=1}^\infty  \frac{\partial}{\partial \alpha}\gamma_k(1,\alpha;v,w)(s-2)^k.
\end{align*}
Then we have
\begin{align*}
    \zeta_2(s, \alpha; v, w) 
    &= 
    \frac{1}{1-s}
    \left\{
    \frac{ \frac{\partial}{\partial \alpha}\gamma_{-1}(1,\alpha;v,w)}{s-2} +  \frac{\partial}{\partial \alpha}\gamma_0(1,\alpha;v,w)
    \right. \\
    &\qquad\qquad\left. +
    \sum_{k=1}^\infty  \frac{\partial}{\partial \alpha}\gamma_k(1,\alpha;v,w)(s-2)^k
    \right\} \\
    &= \left\{\sum_{p=0}^\infty (-1)^{p+1}(s-2)^p \right\} \left\{\sum_{q=-1}^\infty  \frac{\partial}{\partial \alpha}\gamma_q(1,\alpha;v,w)(s-2)^q  \right\}   \\
    &= \frac{1}{vw}\frac{1}{s-2} + \sum_{k=0}^\infty \left\{ \sum_{l=-1}^k (-1)^{k-l+1}  \frac{\partial}{\partial \alpha}\gamma_l(1,\alpha;v,w) \right\}(s-2)^k.
\end{align*}
By the uniqueness of the Laurent series expansion, comparing the coefficients of the $k$-th term, we have
\begin{align*}
    \sum_{l=-1}^k (-1)^{k-l+1}  \frac{\partial}{\partial \alpha}\gamma_l(1,\alpha;v,w) = \gamma_k(2,\alpha;v,w).
\end{align*}
This proves \eqref{th:Main_Theorem2-2} and completes the proof.
\end{proof}

\section{Asymptotic behavior of the Laurent coefficients}

Next, in a similar spirit to \eqref{estim_gamma_k} and
\eqref{estim_gamma_k(a)}, we investigate the asymptotic behavior of
$\gamma_k(2,\alpha;v,w)$ and $\gamma_k(1,\alpha;v,w)$ as
$k\to\infty$.

Define
\begin{equation}\label{H(s)}
H(s)
:=
\zeta_2(s,\alpha;v,w)
-\frac{1}{vw}\frac{1}{s-2}
-\frac{v+w-2\alpha}{2vw}\frac{1}{s-1}.
\end{equation}
Since $\zeta_2(s,\alpha;v,w)$ has no poles other than $s=1$ and
$s=2$, the function $H(s)$ is entire.

For $R>0$, put
\[
M_R^{(2)}
:=
\max_{|s-2|=R}|H(s)|,
\qquad
M_R^{(1)}
:=
\max_{|s-1|=R}|H(s)|.
\]

The following results describe the asymptotic behavior of the
Laurent coefficients at $s=2$ and $s=1$, respectively.

\begin{theorem}\label{th:asymp_s2}
Let $\alpha,v,w>0$. Then, for every $R>0$ and
$k=0,1,2,\ldots$,
\[
\left|
\gamma_k(2,\alpha;v,w)
-
(-1)^k\frac{v+w-2\alpha}{2vw}
\right|
\le
\frac{M_R^{(2)}}{R^k}.
\]
In particular, for every fixed $R>1$,
\[
\gamma_k(2,\alpha;v,w)
=
(-1)^k\frac{v+w-2\alpha}{2vw}
+
O(R^{-k})
\qquad (k\to\infty).
\]
\end{theorem}

\begin{proof}
Recall that $H(s)$ defined by \eqref{H(s)} is entire. Expanding $(s-1)^{-1}$ around $s=2$, we have
\[
\frac{1}{s-1}
=
\sum_{j=0}^{\infty}(-1)^j(s-2)^j
\qquad (|s-2|<1).
\]
Hence, in a neighborhood of $s=2$,
\[
\zeta_2(s,\alpha;v,w)
=
\frac{1}{vw}\frac{1}{s-2}
+
\sum_{j=0}^{\infty}
\left\{
(-1)^j\frac{v+w-2\alpha}{2vw}
+
\frac{H^{(j)}(2)}{j!}
\right\}(s-2)^j.
\]
Comparing this with \eqref{Laurent_s=2}, we obtain
\[
\gamma_k(2,\alpha;v,w)
=
(-1)^k\frac{v+w-2\alpha}{2vw}
+
\frac{H^{(k)}(2)}{k!}.
\]
Since $H(s)$ is entire, Cauchy's integral formula for derivatives gives
\[
\frac{H^{(k)}(2)}{k!}
=
\frac{1}{2\pi i}
\int_{|s-2|=R}\frac{H(s)}{(s-2)^{k+1}}\,ds,
\]
and hence
\[
\left|\frac{H^{(k)}(2)}{k!}\right|
\le
\frac{1}{R^k}\max_{|s-2|=R}|H(s)|
=
\frac{M_R^{(2)}}{R^k}.
\]
This proves the theorem.
\end{proof}

\begin{theorem}\label{th:asymp_s1}
Let $\alpha,v,w>0$. Then, for every $R>0$ and
$k=0,1,2,\ldots$,
\[
\left|
\gamma_k(1,\alpha;v,w)
+\frac{1}{vw}
\right|
\le
\frac{M_R^{(1)}}{R^k}.
\]
In particular, for every fixed $R>1$,
\[
\gamma_k(1,\alpha;v,w)
=
-\frac{1}{vw}
+
O(R^{-k})
\qquad (k\to\infty).
\]
Moreover, the function
\[
\zeta_2(s,\alpha;v,w)
-\frac{v+w-2\alpha}{2vw}\frac{1}{s-1}
\]
is holomorphic on $\{|s-1|<1\}$, and its Taylor expansion at
$s=1$ has radius of convergence exactly $1$.
\end{theorem}

\begin{proof}
The function $H(s)$ is entire, and
\[
\zeta_2(s,\alpha;v,w)
=
\frac{1}{vw}\frac{1}{s-2}
+
\frac{v+w-2\alpha}{2vw}\frac{1}{s-1}
+H(s).
\]
Therefore,
\[
\zeta_2(s,\alpha;v,w)
-\frac{v+w-2\alpha}{2vw}\frac{1}{s-1}
\]
is holomorphic on $\{|s-1|<1\}$. Since $s=2$ lies on the boundary of this disc and is a pole, its Taylor expansion at $s=1$ has radius of convergence exactly $1$.
For $|s-1|<1$, we have
\[
\frac{1}{s-2}
=-\sum_{j=0}^{\infty}(s-1)^j.
\]
Expanding $H(s)$ at $s=1$ and comparing with \eqref{Laurent_s=1}, we obtain
\[
\gamma_k(1,\alpha;v,w)
=
-\frac{1}{vw}
+\frac{H^{(k)}(1)}{k!}.
\]
By Cauchy's integral formula,
\[
\left|\frac{H^{(k)}(1)}{k!}\right|
\le
\frac{1}{R^k}\max_{|s-1|=R}|H(s)|
=
\frac{M_R^{(1)}}{R^k}.
\]
This proves the theorem.
\end{proof}

\begin{remark}
From the above two theorems, we obtain
\[
\lim_{k\to\infty}(-1)^k \gamma_k(2,\alpha;v,w)=\frac{v+w-2\alpha}{2vw}, \qquad
\lim_{k\to\infty}\gamma_k(1,\alpha;v,w)
=
-\frac{1}{vw}.
\]
In contrast to \eqref{estim_gamma_k} and \eqref{estim_gamma_k(a)}, both sequences converge to simple closed-form expressions.
\end{remark}

\section{A Mellin-type representation}

In this section, we derive integral representations for the
Laurent coefficients at $s=2$ using the Mellin transform.
Recall the Mellin integral representation
\[
\zeta_2(s,\alpha;v,w)
=
\frac{1}{\Gamma(s)}
\int_0^\infty
\frac{t^{s-1}e^{-\alpha t}}
{(1-e^{-vt})(1-e^{-wt})}\,dt,
\qquad \Re(s)>2.
\]
Define
\[
R_{\alpha,v,w}(t)
:=
\frac{e^{-\alpha t}}
{(1-e^{-vt})(1-e^{-wt})}
-
\frac{e^{-t}}{vw\,t^2}.
\]
Since
\[
R_{\alpha,v,w}(t)=O(t^{-1})
\qquad (t\to0^+),
\]
and $R_{\alpha,v,w}(t)$ decays exponentially as $t\to\infty$,
the integral
\[
\int_0^\infty
t^{s-1}R_{\alpha,v,w}(t)\,dt
\]
is holomorphic in a neighborhood of $s=2$.
\begin{proposition}\label{prop:integral-gamma-s2}
For every integer $k\geq0$,
\begin{align}
\gamma_k(2,\alpha;v,w)
&=
\frac{(-1)^{k+1}}{vw}
+
\frac{1}{k!}
\left.
\frac{d^k}{ds^k}
\left\{
\frac{1}{\Gamma(s)}
\int_0^\infty
t^{s-1}R_{\alpha,v,w}(t)\,dt
\right\}
\right|_{s=2}.
\label{eq:integral-gamma-s2}
\end{align}
Equivalently,
\begin{align}
\gamma_k(2,\alpha;v,w)
&=
\frac{(-1)^{k+1}}{vw}
+
\sum_{j=0}^k
\frac{1}{j!(k-j)!}
\left.
\frac{d^{\,k-j}}{ds^{\,k-j}}
\frac{1}{\Gamma(s)}
\right|_{s=2}
\nonumber\\
&\qquad\qquad\times
\int_0^\infty
t\,R_{\alpha,v,w}(t)
(\log t)^j\,dt.
\label{eq:integral-gamma-s2-explicit}
\end{align}
\end{proposition}
\begin{proof}
By the definition of $R_{\alpha,v,w}(t)$, we have
\begin{align*}
\zeta_2(s,\alpha;v,w)
&=
\frac{1}{\Gamma(s)}
\int_0^\infty
t^{s-1}R_{\alpha,v,w}(t)\,dt
+
\frac{1}{vw\,\Gamma(s)}
\int_0^\infty
t^{s-3}e^{-t}\,dt
\\
&=
\frac{1}{\Gamma(s)}
\int_0^\infty
t^{s-1}R_{\alpha,v,w}(t)\,dt
+
\frac{\Gamma(s-2)}{vw\,\Gamma(s)}
\\
&=
\frac{1}{\Gamma(s)}
\int_0^\infty
t^{s-1}R_{\alpha,v,w}(t)\,dt
+
\frac{1}{vw(s-1)(s-2)}.
\end{align*}
Putting $s=2+\varepsilon$, we obtain
\[
\frac{1}{(s-1)(s-2)}
=
\frac{1}{\varepsilon(1+\varepsilon)}
=
\frac{1}{\varepsilon}
+
\sum_{k=0}^{\infty}(-1)^{k+1}\varepsilon^k.
\]
Comparing this with the Laurent expansion
\[
\zeta_2(s,\alpha;v,w)
=
\frac{1}{vw}\frac{1}{s-2}
+
\sum_{k=0}^{\infty}
\gamma_k(2,\alpha;v,w)(s-2)^k,
\]
we obtain \eqref{eq:integral-gamma-s2}.
Finally, applying Leibniz's rule and using
\[
\left.
\frac{d^j}{ds^j}t^{s-1}
\right|_{s=2}
=
t(\log t)^j,
\]
gives \eqref{eq:integral-gamma-s2-explicit}.
\end{proof}
In particular, taking $k=0$ in Proposition
\ref{prop:integral-gamma-s2}, we obtain the representation \eqref{eq:integral-gamma0-s2}
for the constant term.
\begin{corollary}\label{cor:integral-gamma0-s2}
We have
\begin{align}
\gamma_0(2,\alpha;v,w)
&=
-\frac{1}{vw}
+
\int_0^\infty
\left\{
\frac{t e^{-\alpha t}}
{(1-e^{-vt})(1-e^{-wt})}
-
\frac{e^{-t}}{vw\,t}
\right\}\,dt.
\label{eq:integral-gamma0-s2}
\end{align}
\end{corollary}

\section{Special values and arithmetic aspects}

We conclude with some simple special values of the constant terms at
$s=1$ and $s=2$.
These examples are obtained by reducing the Barnes double zeta function
to finite linear combinations of Hurwitz zeta functions.

As a special case of the reduction formula of Sakane and Aoki \cite{SakaneAoki2022}*{Theorem 4} for the Barnes multiple zeta function with positive rational periods, we obtain the following expression. Taking
$
N=2,\qquad x=\alpha,\qquad w_1=1,\qquad w_2=q,
$
where $q$ is a positive integer, their formula reduces to
\begin{align}
\zeta_2(s,\alpha;1,q)
&=
q^{-s}\sum_{r=0}^{q-1}
\zeta_H\left(s-1,\frac{\alpha+r}{q}\right)
\nonumber\\
&\quad+
q^{-s-1}\sum_{r=0}^{q-1}
(q-r-\alpha)
\zeta_H\left(s,\frac{\alpha+r}{q}\right).
\label{eq:reduction-rational-ratio}
\end{align}
Thus $q^s\zeta_2(s,\alpha;1,q)$ is a finite linear combination of
Hurwitz zeta functions with coefficients independent of $s$.

We first consider the constant term at $s=1$.
With $v=1$ and $w=q$, a particularly simple family is obtained by taking
\[
\alpha=\frac{q+1}{2}.
\]
In this case, the residue at $s=1$ vanishes.

\begin{proposition}\label{prop:s1-special-values}
Let $q\geq2$ be an integer. Then
\begin{align}
\gamma_0\left(1,\frac{q+1}{2};1,q\right)
=
-\frac12
+
\frac{1}{q}\left\lfloor\frac q2\right\rfloor
-
\frac{\pi}{2q^2}
\sum_{r=0}^{\lfloor q/2\rfloor-1}
(q-1-2r)
\tan\left(\frac{(2r+1)\pi}{2q}\right).
\label{eq:gamma0-s1-q}
\end{align}
\end{proposition}

\begin{proof}
Putting $\alpha=(q+1)/2$ in
\eqref{eq:reduction-rational-ratio}, we obtain
\begin{align*}
\zeta_2\left(s,\frac{q+1}{2};1,q\right)
&=
q^{-s}
\sum_{r=0}^{q-1}
\zeta_H\left(
s-1,\frac{q+1+2r}{2q}
\right)
\\
&\quad+
q^{-s-1}
\sum_{r=0}^{q-1}
\frac{q-1-2r}{2}
\zeta_H\left(
s,\frac{q+1+2r}{2q}
\right).
\end{align*}
Since
\[
\sum_{r=0}^{q-1}(q-1-2r)=0,
\]
the pole at $s=1$ in the second sum cancels. Hence, by taking the
constant term at $s=1$, we have
\begin{align}
\gamma_0\left(1,\frac{q+1}{2};1,q\right)
&=
\frac1q
\sum_{r=0}^{q-1}
\zeta_H\left(
0,\frac{q+1+2r}{2q}
\right)
\nonumber\\
&\quad+
\frac1{q^2}
\sum_{r=0}^{q-1}
\frac{q-1-2r}{2}
\gamma_0\left(
\frac{q+1+2r}{2q}
\right).
\label{eq:gamma0-s1-before-pairing}
\end{align}
Using
\[
\zeta_H(0,a)=\frac12-a,
\]
the first sum on the right-hand side is
\[
\frac1q
\sum_{r=0}^{q-1}
\left(
\frac12-\frac{q+1+2r}{2q}
\right)
=
-\frac12.
\]
We now consider the second sum. Put
$
a_r=(q+1+2r)/(2q).
$
For $0\leq r\leq\lfloor q/2\rfloor-1$, pairing the terms corresponding
to $r$ and $q-1-r$, we have
$
a_{q-1-r}=2-a_r
$
and
$
q-1-2(q-1-r)=-(q-1-2r).
$
Since
\[
\gamma_0(a)=-\psi(a),
\]
where $\psi=\Gamma'/\Gamma$ is the digamma function, the recurrence and reflection formulas
\[
\psi(z+1)=\psi(z)+\frac1z,
\qquad
\psi(1-z)-\psi(z)=\pi\cot(\pi z)
\]
give
\begin{align*}
\gamma_0(a)-\gamma_0(2-a)
&=
\psi(2-a)-\psi(a)
=
\frac1{1-a}+\pi\cot(\pi a).
\end{align*}
For $a=a_r$, we have
\[
1-a_r=\frac{q-1-2r}{2q},
\]
and therefore
\[
\frac{q-1-2r}{2}
\frac1{1-a_r}
=q.
\]
Moreover,
\[
a_r
=
\frac12+\frac{2r+1}{2q},
\]
so that
\[
\cot(\pi a_r)
=
-\tan\left(\frac{(2r+1)\pi}{2q}\right).
\]
Consequently, each pair contributes
\[
q
-
\frac{\pi}{2}(q-1-2r)
\tan\left(\frac{(2r+1)\pi}{2q}\right).
\]
There are $\lfloor q/2\rfloor$ such pairs. When $q$ is odd,
the remaining middle term has coefficient zero and hence does not
contribute. Substituting these expressions into
\eqref{eq:gamma0-s1-before-pairing}, we obtain
\eqref{eq:gamma0-s1-q}.
\end{proof}

\medskip

\begin{theorem}\label{th:transcendence-s1}
For every integer $q\geq2$,
\[
\gamma_0\left(1,\frac{q+1}{2};1,q\right)
\]
is transcendental. In particular, it is irrational.
\end{theorem}

\begin{proof}
For
\[
0\leq r\leq
\left\lfloor\frac q2\right\rfloor-1,
\]
we have
\[
0<
\frac{(2r+1)\pi}{2q}
<
\frac{\pi}{2},
\]
and hence
\[
(q-1-2r)
\tan\left(\frac{(2r+1)\pi}{2q}\right)>0.
\]
Therefore
\[
\sum_{r=0}^{\lfloor q/2\rfloor-1}
(q-1-2r)
\tan\left(\frac{(2r+1)\pi}{2q}\right)
\neq0.
\]
Furthermore, each number
\[
\tan\left(\frac{(2r+1)\pi}{2q}\right)
\]
is algebraic, since it can be expressed in terms of roots of unity.
Thus the finite sum above is a nonzero algebraic number.
It follows from \eqref{eq:gamma0-s1-q} that
\[
\gamma_0\left(1,\frac{q+1}{2};1,q\right)
=
A+B\pi,
\]
where $A$ is rational and $B$ is a nonzero algebraic number.
Since $B\neq0$, algebraicity of $A+B\pi$ would imply algebraicity of
$\pi$, contradicting its transcendence. Thus $A+B\pi$ is transcendental.
This proves the assertion.
\end{proof}

\begin{example}
For $q=2$, Proposition \ref{prop:s1-special-values} gives
\[
\gamma_0\left(1,\frac32;1,2\right)
=
-\frac{\pi}{8}.
\]
For $q=3$, it gives
\[
\gamma_0(1,2;1,3)
=
-\frac16-\frac{\pi}{9\sqrt3}.
\]
These are the cases $q=2,3$ of Theorem \ref{th:transcendence-s1}.
Their transcendence also follows directly from these explicit expressions,
since $\pi$ is transcendental and the coefficients of $\pi$ are nonzero
algebraic numbers.
\end{example}

We next consider the constant term at $s=2$.

\begin{example}
For $\alpha=v=w=1$, we have
\[
\zeta_2(s,1;1,1)=\zeta(s-1).
\]
Hence
\[
\gamma_0(1,1;1,1)=\zeta(0)=-\frac12,
\]
while the constant term at $s=2$ is
\[
\gamma_0(2,1;1,1)=\gamma.
\]
\end{example}

When $v=w=1$, the Barnes double zeta function can be written in a particularly
simple form:
\begin{align}
\zeta_2(s,\alpha;1,1)
&=
\sum_{m=0}^{\infty}\sum_{n=0}^{\infty}
\frac{1}{(\alpha+m+n)^s}
\nonumber\\
&=
\sum_{k=0}^{\infty}
\frac{k+1}{(\alpha+k)^s}
\nonumber\\
&=
\zeta_H(s-1,\alpha)
+
(1-\alpha)\zeta_H(s,\alpha).
\label{eq:zeta2-hurwitz-vw1}
\end{align}

\begin{proposition}\label{prop:s2-integer-values}
For every positive integer $N$, we have
\begin{align}
\gamma_0(2,N;1,1)
&=
\gamma
-
\sum_{m=1}^{N-1}\frac1m
+
(N-1)\sum_{m=1}^{N-1}\frac1{m^2}
-\frac{N-1}{6}\pi^2.
\label{eq:gamma0-s2-integer}
\end{align}
\end{proposition}

\begin{proof}
From \eqref{eq:zeta2-hurwitz-vw1}, we have
\[
\zeta_2(s,\alpha;1,1)
=
\zeta_H(s-1,\alpha)
+
(1-\alpha)\zeta_H(s,\alpha).
\]
Since
\[
\zeta_H(s-1,\alpha)
=
\frac{1}{s-2}
+
\gamma_0(\alpha)
+
O(s-2)
\]
near $s=2$, comparison with the Laurent expansion of
$\zeta_2(s,\alpha;1,1)$ at $s=2$ gives
\begin{align}
\gamma_0(2,\alpha;1,1)
=
\gamma_0(\alpha)
+
(1-\alpha)\zeta_H(2,\alpha).
\label{eq:gamma0-s2-hurwitz}
\end{align}

For a positive integer $N$, the shift relation for the Hurwitz zeta function gives
\[
\zeta_H(s,N)
=
\zeta(s)
-
\sum_{m=1}^{N-1}\frac1{m^s}.
\]
Taking the constant term at $s=1$, we obtain
\[
\gamma_0(N)
=
\gamma
-
\sum_{m=1}^{N-1}\frac1m.
\]
Also, putting $s=2$, we obtain
\[
\zeta_H(2,N)
=
\frac{\pi^2}{6}
-
\sum_{m=1}^{N-1}\frac1{m^2}.
\]
Substituting these formulas into \eqref{eq:gamma0-s2-hurwitz}, we obtain
\eqref{eq:gamma0-s2-integer}.
\end{proof}

\begin{example}
The first few cases of Proposition \ref{prop:s2-integer-values} are
\[
\gamma_0(2,1;1,1)=\gamma,
\]
\[
\gamma_0(2,2;1,1)
=
\gamma-\frac{\pi^2}{6},
\]
and
\[
\gamma_0(2,3;1,1)
=
\gamma+1-\frac{\pi^2}{3}.
\]
Thus, although the arithmetic nature of the Euler constant $\gamma$ itself
is still unknown, differences between distinct members of this family can
have a completely determined arithmetic nature.
\end{example}

\begin{theorem}\label{th:transcendence-s2}
Among the numbers
\[
\gamma_0(2,N;1,1)
\qquad (N=1,2,3,\ldots),
\]
all but at most one are transcendental.
\end{theorem}

\begin{proof}
Let $M$ and $N$ be distinct positive integers. By
\eqref{eq:gamma0-s2-integer},
\begin{align*}
&\gamma_0(2,N;1,1)-\gamma_0(2,M;1,1)
\\
&\quad=
-\sum_{m=1}^{N-1}\frac1m
+(N-1)\sum_{m=1}^{N-1}\frac1{m^2}
+\sum_{m=1}^{M-1}\frac1m
-(M-1)\sum_{m=1}^{M-1}\frac1{m^2}
-\frac{N-M}{6}\pi^2.
\end{align*}
The terms other than the last one are rational, while
$(N-M)/6$ is a nonzero rational number. Since $\pi^2$ is transcendental,
the difference
$
\gamma_0(2,N;1,1)-\gamma_0(2,M;1,1)
$
is transcendental.
If both $\gamma_0(2,N;1,1)$ and $\gamma_0(2,M;1,1)$ were algebraic,
then their difference would also be algebraic, which is impossible.
Hence at most one member of the family is algebraic.
\end{proof}

\begin{remark}
The above results exhibit two different arithmetic phenomena.
At $s=1$, the values in Theorem \ref{th:transcendence-s1} are rational
numbers plus nonzero algebraic multiples of $\pi$, and hence every member
of this family is transcendental.
At $s=2$, the values $\gamma_0(2,N;1,1)$ involve the Euler constant
$\gamma$, whose irrationality is still unknown; nevertheless,
Proposition \ref{prop:s2-integer-values} shows that the differences between
distinct members of this family are transcendental. Consequently,
all but at most one member of this family are transcendental, as stated
in Theorem \ref{th:transcendence-s2}.

These examples are not intended to provide a general irrationality or
transcendence criterion for the Laurent coefficients. Formula
\eqref{eq:reduction-rational-ratio} and the special case
\eqref{eq:zeta2-hurwitz-vw1} suggest that further arithmetic questions may
be studied when $v/w$ is rational.
\end{remark}

\section*{Acknowledgments}
At the 2025 Annual Meeting of the Mathematical Society of Japan, held at Waseda University, I received valuable comments from Prof. Genki Shibukawa, which contributed to the further development of this paper. I would like to take this opportunity to express my deep gratitude.
I am also deeply grateful to Prof. Akihiko Yukie, Prof. Masatoshi Suzuki, and Prof. Hirotaka Akatsuka for their valuable insights and advice.
I would like to thank all members of my laboratory and all members of the Kansai Multiple Zeta Study Group.

\bibliographystyle{amsalpha}
\bibliography{References} 

\end{document}